\documentclass{article}

\usepackage[vmargin=3cm,hmargin=3cm,headheight=18pt,top=3cm,headsep=.5cm]{geometry}

\usepackage{amssymb}
\usepackage{latexsym}
\usepackage{epsfig}
\usepackage{graphicx} % Permite incluir dibujos jpg,bmp,eps...
\usepackage{color}
\usepackage{amsmath}
\usepackage{amsfonts}
\usepackage{url,hyperref}

\usepackage{tikz,pgfplots}
\usepackage{tikz-cd}
\usepackage{stackrel}

\newtheorem{theorem}{Theorem}
\newtheorem{lemma}{Lemma}
\newtheorem{remark}{Remark}
\newtheorem{corollary}{Corollary}

\pgfplotsset{compat=1.16} 

\begin{document}

\title{Using Aichinger's equation to characterize polynomial functions}

%\title{The double sidedness of the matrix inverse: a proof by induction}
%\title{$A\cdot B=I \Longleftrightarrow I=B\cdot A$ }

\author{J. M. Almira\\
\normalsize {Departamento de Ingenier\'{\i}a y Tecnolog\'{\i}a de Computadores,}\\
\normalsize { \'Area de Matem\'atica Aplicada, Facultad de Inform\'atica,}\\
\normalsize{Universidad de Murcia, Campus Espinardo, 30100 Murcia, Spain.}\\
\normalsize {e-mail: jmalmira@um.es}\\
}
 %%leave blank in initial submission to allow for double blind reviewing

\date{}
\maketitle

%\begin{center} To Professors Maciej Sablik and László Székelyhidi, in their 70th birthady, with admiration and friendship. \end{center}

\begin{abstract} Aichinger's equation is used to give simple proofs of several well-known characterizations of polynomial functions as solutions of certain functional equations. Concretely, we use that Aichinger's equation characterizes polynomial functions to solve, for arbitrary commutative groups, Ghurye-Olkin's functional equation, Wilson's functional equation, the Kakutani-Nagumo-Walsh functional equation, and a general version of Fréchet's unmixed functional equation.

\end{abstract}

\medbreak 

%\noindent \thanks{  {\it 2020 MSC:} ??} \medbreak

\noindent   \vspace*{0.5cm}

\section{Introduction}

%
%We consider Aichinger's equation $$f(x_1+\cdots+x_{m+1})=\sum_{i=1}^{m+1}g_i(x_1,x_2,\cdots, \widehat{x_i},\cdots, x_{m+1})$$ for functions defined on commutative semigroups which take values on commutative groups. The solutions of this equation are, under very mild hypotheses, generalized polynomials. 

Given a commutative semigroup $(S,+)$ and a commutative group $(H,+)$, we say that $f:S\to H$ is a polynomial function (also named generalized polynomial) of degree $\leq m$ if $f$ solves Fréchet's mixed differences functional equation:
$$
\Delta_{h_1}\Delta_{h_2}\cdots \Delta_{h_{m+1}}f(x)=0\text{ for all } h_1,\cdots, h_{m+1},x\in S.
$$
When we evaluate a polynomial function of degree $\leq m$ on a sum $u_1+u_2+\cdots+u_{m+1}$ of $m+1$ variables $u_1,\cdots, u_{m+1}$, we obtain a sum of functions with the property that each one of them depends on at most $m$ variables $u_{i_1},\cdots,u_{i_m}$ ($i_k\in\{1,\cdots,m+1\}$ for all $k$). Indeed, this fact completely characterises polynomial functions, as was proved by Aichinger and Moosbauer \cite{AM} for functions defined on abelian groups and by Almira \cite{almira} for functions defined on commutative semigroups $S$ which satisfy that $S+S=S$ and $0\in S$. Concretely \cite[Lemma 4.1]{AM} states that, when $(S,+)$ is also a group, the function $f:S\to H$ is a generalized polynomial of degree $\leq m$ if and only if it solves Aichinger's equation: 
\begin{equation}\label{aichinger}
f(x_1+\cdots+x_{m+1})=\sum_{i=1}^{m+1}g_i(x_1,x_2,\cdots, \widehat{x_i},\cdots, x_{m+1})
\end{equation}
for certain functions $g_i:S^{m}\to H$, $i=1,2,\cdots, m+1$. Here $\widehat{x_{i}}$ means that $g_i$  does not depend on $x_i$. Moreover, \cite[Theorem 2]{almira} proves that, if $S+S=S$ and $f$ solves \eqref{aichinger}, then $f$ is a polynomial function of degree $\leq m$. Moreover, if $S$ also satisfies that $0\in S$, then all polynomial functions of degree $\leq m$ are solutions of \eqref{aichinger}. 

The main goal of this note is to use these results to give simple proofs of several well-known (and useful) characterizations of polynomial functions as solutions of certain functional equations. Concretely, we use that Aichinger's equation characterizes polynomial functions to solve Ghurye-Olkin's functional equation (Corollary \ref{GO}), Wilson's functional equation (Corollary \ref{Wilson}), the Kakutani-Nagumo-Walsh functional equation (Corollary \ref{KNW}), and a general version of Fréchet's unmixed functional equation (Corollary \ref{GFFE}). In all cases we give proofs for general commutative groups or semigroups and we do not worry about the regularity properties of polynomial functions. Although the results we prove are well-known, our proofs are surprisingly easier than the original ones.

\section{Characterizations of polynomial functions}
We use Aichinger's equation to prove the main result of this note, which is the following theorem:
\begin{theorem} \label{main}
Assume that  $(S,+)$ is a commutative semigroup such that $S+S=S$, $(R,+,\cdot)$ is a commutative ring, $c:S\to S$ is an automorphism, and $f:(S,+)\to (R,+)$ is a map that satisfies 
\begin{equation}\label{GOcero}
f(x+c(y))=\sum_{j=1}^Np_j(x)a_j(y)+\sum_{k=1}^Mq_k(y)b_k(x) 
\end{equation}
where $p_j,q_k:(S,+)\to (R,+)$ are polynomial functions and $\deg(p_j)\leq r$, $\deg(q_k)\leq s$ for all $j,k$. Then $f:(S,+)\to (R,+)$ is a polynomial function and $\deg (f) \leq r+s+1$. 

\end{theorem}
This theorem is a generalized version of \cite[Lemma 2.1]{almirajmaa} and is the key to proving Corollary \ref{GO}, which is a generalization of a result proved by Ghurye and Olkin in \cite[Lemma 3]{Gu_O} (see also \cite[Theorem 1.3]{almirajmaa}), where the equation was used for the characterization of Gaussian probability distributions (see also \cite[Chapter 7]{MaPe}). Our proof  is not based on the very technical tools used in \cite{Gu_O}. Moreover, it simplifies the proof given in \cite{almirajmaa}, since it avoids using  exponential polynomials. In particular, it does not depend on \cite[Lemma 4.3]{Sz1} and/or \cite[Theorem 5.10]{Sz}.

\noindent \textbf{Proof of Theorem \ref{main}.}  
Take $x_1,\cdots ,x_{r+1},y_1,\cdots,y_{s+1}$, a set of $r+s+2$ variables, then:
\begin{eqnarray*}
& \ & f(x_1+\cdots+x_{r+1}+y_1+\cdots+y_{s+1}) \\ 
&=&  f(x_1+\cdots+x_{r+1}+c(c^{-1}(y_1)+\cdots+c^{-1}(y_{s+1})))\\
&=& \sum_{j=1}^Np_j(x_1+\cdots+x_{r+1})a_j(c^{-1}(y_1)+\cdots+c^{-1}(y_{s+1}))\\ &\ & \ +\sum_{k=1}^Mq_k(c^{-1}(y_1)+\cdots+c^{-1}(y_{s+1}))b_k(x_1+\cdots+x_{r+1}) ,
\end{eqnarray*}
now, each term $p_j(x_1+\cdots+x_{r+1})$ can be decomposed as a sum of terms, each one depending on at most $r$ of the variables  $x_1,\cdots,x_{r+1}$; and each term $q_k(c^{-1}(y_1)+\cdots+c^{-1}(y_{s+1}))$ can be decomposed as a sum of terms, each one depending on at most $s$ of the variables  $y_1,\cdots,y_{s+1}$. Thus, $f$ satisfies Aichinger's equation of order $r+s+1$, which means that $f$ is a polynomial function of degree $\leq r+s+1$. 

{\hfill $\Box$}

Aichinger's equation can also be used to prove the following 

\begin{lemma}\label{lemma}
Assume that  $(S,+)$ is a commutative semigroup such that $0\in S= S+S$, $(H,+)$ is a commutative group, $c:S\to S$ is an automorphism, and $f:(S,+)\to (H,+)$ is a map such that $g(x)=f(c(x))$ is a polynomial function of degree $\leq m$. Then $f$ is also a polynomial function of degree $\leq m$.
\end{lemma}

\noindent \textbf{Proof. }
Given a set of $m+1$ variables $x_1,\cdots, x_{m+1}$, we have that 
\begin{eqnarray*}
f(x_1+\cdots+x_{m+1})&=& f(c(c^{-1}(x_1)+\cdots+c^{-1}(x_{m+1})))\\
&=&g(c^{-1}(x_1)+\cdots+c^{-1}(x_{m+1}))\\
&=&\sum_{i=1}^{m+1}g_i( c^{-1}(x_1),\cdots,\widehat{c^{-1}(x_i)},\cdots, c^{-1}(x_{m+1}))\\
&=&\sum_{i=1}^{m+1}G_i( x_1,\cdots,\widehat{x_i},\cdots x_{m+1}),
\end{eqnarray*}
which means that $f$ satisfies Aichinger's equation of order $m$.  {\hfill $\Box$}

We are now in a good position to study Ghurye-Olkin's functional equation:

\begin{corollary}[Ghurye-Olkin's functional equation]\label{GO}
With the same notation as used in Theorem \ref{main}, assume that $(S,+)$ is a commutative group, $c_i:S\to S$, $i=1,\cdots,n$ are automorphisms and that $c_i-c_j$ is also an automorphism whenever $i\neq j$. Let $f_i:S\to R$, $i=1,\cdots,n$  be such that 
\begin{equation}\label{eqGO}
\sum_{i=1}^nf_i(x+c_i(y))=\sum_{j=1}^Np_j(x)a_j(y)+\sum_{k=1}^Mq_k(y)b_k(x) 
\end{equation}
where $p_j,q_k:(S,+)\to (R,+)$ are polynomial functions and $\deg(p_j)\leq r$, $\deg(q_k)\leq s$ for all $j,k$. Then each function $f_i:(S,+)\to (R,+)$ is a polynomial function and $\deg (f_i) \leq r+s+n$, $i=1,\cdots,n$. Moreover, in the particular case that $p_j=q_k=0$ for all $j,k$ (so that the second member of the equation is $0$), the equation is well defined for functions $f_i$ taking values on a commutative group and each $f_i$ that solves the equation is a polynomial function of degree at most $n-1$. 
\end{corollary}
\noindent \textbf{Proof. } Assume that the second member of \eqref{eqGO} is not $0$. Theorem \ref{main} is case $n=1$, so that it is natural to proceed by induction on $n$. Indeed, we prove by induction on $n$ that $f_n$ is a polynomial function of degree $\leq r+s+n$, and a simple rearrangement of the functions $f_1,\cdots,f_n$ proves the result for all functions $f_i$, $i=1,\cdots,n$.  

Let us assume that $n>1$ and the result holds for $n-1$. If we consider both members of the equation as functions $F(x,y)$ defined on $S\times S$, we can apply the difference operator $$\Delta_{(h_1,-c_1^{-1}(h_1))}F(x,y)=F(x+h_1,y-c_1^{-1}(h_1))-F(x,y)$$ to both sides of the equation to get that 
$$
\sum_{i=2}^ng_i(x+c_i(y))=\sum_{j=1}^Np_j^*(x)a_j^*(y)+\sum_{k=1}^Mq_k^*(y)b_k^*(x), 
$$
where
\begin{eqnarray*}
g_i(x+c_i(y)) &=& \Delta_{(h_1,-c_1^{-1}(h_1))}f_i(x+c_i(y)) \\
&=& f_i(x+h_1+c_i(y-c_1^{-1}(h_1))) - f_i(x+c_i(y)) \\
&=& f_i(x+h_1+c_i(y)-c_i(c_1^{-1}(h_1)))-f_i(x+c_i(y)) \\
&=& \Delta_{h_1-c_i(c_1^{-1}(h_1))}f_i(x+c_i(y))\\
&=& \Delta_{(1_d-c_i\circ c_1^{-1})(h_1)}f_i(x+c_i(y))
\end{eqnarray*}
and $p_j^*,q_k^*:(S,+)\to (R,+)$ are polynomial functions satisfying $\deg(p_j^*)\leq r$, $\deg(q_k^*)\leq s$ for all $j,k$. Thus, the induction hypothesis implies that all functions $g_i$, $i=2,\cdots, n$ are polynomial functions of degree $\leq r+s+n-1$. In particular, for each $h_1\in S$, $\Delta_{(1_d-c_i\circ c_1^{-1})(h_1)}f_n$ is a polynomial function of degree $\leq r+s+n-1$. Now, $(c_1-c_i)\circ c_1^{-1}= 1_d-c_i\circ c_1^{-1}$ is bijective, so that, for each $h\in S$, $\Delta_{h}f_n$ is a polynomial function of degree $\leq r+s+n-1$. In particular, $f_n$ satisfies Fréchet's mixed functional equation 
$$
\Delta_{h_1}\Delta_{h_2}\cdots \Delta_{h_{r+s+n}} \Delta_{h}f(x)=0\text{ for all } h_1,\cdots, h_{r+s+n-1},h,x\in S.
$$
Hence $f_n$ is a polynomial function of degree $\leq r+s+n$.

Let us now assume that that $p_j=q_k=0$ for all $j,k$ and that all functions $f_i$ take values on a commutative group (not necessarily a ring). Then the equation is of the form 
\begin{equation}\label{eqGOcase0}
\sum_{i=1}^nf_i(x+c_i(y))=0.
\end{equation}
Again,  we prove by induction on $n$ that $f_n$ is a polynomial function of degree $\leq n-1$, and a simple rearrangement of the functions $f_1,\cdots,f_n$ proves the result for all functions $f_i$, $i=1,\cdots,n$.
For $n=1$, the equation becomes $f_1(x+c_1(y))=0$ and taking $y=0$ we get that $f_1(x)=0$ and $f_1$ is a polynomial of degree $0$. We assume that the result holds true for $n-1$ and consider the case $n>1$. Applying to both sides of \eqref{eqGOcase0} the operator  $\Delta_{(h_1,-c_1^{-1}(h_1))}$ we obtain that 
\[
\sum_{i=2}^n\Delta_{(1_d-c_i\circ c_1^{-1})(h_1)}f_i(x+c_i(y))=0. 
\]
Thus, the induction hypothesis and the fact that $1_d-c_i\circ c_1^{-1}$ is an automorphism for each $i\geq 2$ imply that $\Delta_{h}f_i$ is a polynomial function of degree at most $n-2$ for each $i=2,\cdots,n$ and for every $h\in S$. In particular, $f_n$ is a polynomial function of degree at most $n-1$. 
{\hfill $\Box$}
\begin{remark} \label{rem}
The transformation $\widetilde{f_i}(x)=f_i(\beta_i(x))$ reduces the equation 
\begin{equation}\label{GOM}
\sum_{i=1}^nf_i(\beta_i(x)+\delta_i(y))=\sum_{j=1}^Np_j(x)a_j(y)+\sum_{k=1}^Mq_k(y)b_k(x) 
\end{equation}
to the equation 
\begin{equation}\label{GOM2}
\sum_{i=1}^n\widetilde{f_i}(x+(\beta_i^{-1}\circ \delta_i)(y))=\sum_{j=1}^Np_j(x)a_j(y)+\sum_{k=1}^Mq_k(y)b_k(x), 
\end{equation}
so that, if $\beta_i,\delta_i:S\to S$ are automorphisms such that $\beta_i^{-1}\circ \delta_i-\beta_j^{-1}\circ \delta_j$ is invertible whenever $i\neq j$ and $p_j,q_k:(S,+)\to (R,+)$ are polynomial functions, $\deg(p_j)\leq r$, $\deg(q_k)\leq s$ for all $j,k$, then we can use Corollary \ref{GO} with $c_i=\beta_i^{-1}\circ \delta_i$ and Lemma 1 to conclude that the functions $f_i$ that solve equation \eqref{GOM} are polynomial functions of degree at most $r+s+n$. Moreover, if the second member of  \eqref{GOM} is  $0$, the functions $f_i$ that solve the equation are polynomial functions of degree at most $n-1$. 
\end{remark}
The following equation was studied one hundred years ago by Wilson \cite{wilson}: 
\begin{corollary}[Wilson's functional equation]\label{Wilson}
Assume that $(S,+)$, $(H,+)$  are commutative groups and $\beta_i,\delta_i:S\to S$ are automorphisms such that $\beta_i^{-1}\circ \delta_i-\beta_j^{-1}\circ \delta_j$ is invertible whenever $i\neq j$. Assume also that the functions $f_i,a,b:S\to H$ ($i=1,\cdots,n$) solve the equation
\[
\sum_{i=1}^nf_i(\beta_i(x)+\delta_i(y))=a(x)+b(y)
\]
Then $a,b$ are polynomial functions of degree at most $n$.
\end{corollary}
\noindent \textbf{Proof.} A direct application of Remark \ref{rem} with $r=s=0$ shows that all functions $f_i$, $i=1,\cdots,n$ are polynomial functions of degree at most $n$. Hence $a$, $b$ are also polynomial functions of degree at most $n$. {\hfill $\Box$}

The following equation was studied by Almira in \cite{almirajmaa, A_JJA}.
\begin{corollary}[Generalized Fréchet's unmixed functional equation] \label{GFFE} Let  $(S,+)$ and $(H,+)$  be commutative groups and assume that $f:S\to H$ solves the equation 
\begin{equation}
q(\tau_h)(f)=\sum_{k=0}^na_kf(x+kh)=0 \text{ for all } x,h\in S,
\end{equation} 
where $q(z)=a_0+\cdots+a_nz^n$ is a polynomial and $\tau_h:S\to S$ is the translation operator, $\tau_h(x)=x+h$. Then $f$ is a polynomial function of degree at most $s=\#\{k:a_k\neq 0\}-1$.
\end{corollary}
\noindent  \textbf{Proof. } Set $\{k:a_k\neq 0\}=\{k_1,k_2,\cdots,k_s\}$ and apply Corollary \ref{GO} with $p_j=q_k=0$ for all $j,k$, $f_i=a_{k_i}f(x)$ and $c_i(x)=k_ix:=x+\cdots^{k_i\text{ times }}+x$, $i=1,\cdots,s$.  {\hfill $\Box$}

The following equation was introduced by S. Kakutani and M. Nagumo \cite{kn} and J. L.  Walsh \cite{wal}  in the 1930's. The equation was extensively studied by S. Haruki \cite{h1,h2,h3} in the 1970's  and  1980's.  
\begin{corollary} \label{KNW}
Let $(H,+)$ be a commutative group, let $f:\mathbb{C}\to H$ be a solution of the Kakutami-Nagumo-Walsh functional equation
\begin{equation}\label{eqKNW}
\frac{1}{N}\sum_{k=0}^{N-1}f(z+w^kh)=f(z) \text{ for all } z,h\in\mathbb{C},
\end{equation}
where $w$ is any primitive $N$-th root of $1$. Then $f$ is a polynomial function of degree $\leq N$. 
 \end{corollary}

\noindent  \textbf{Proof. } Use Corollary \ref{Wilson} with $n=N$, $f_i(z)=\frac{1}{N}f(z)$, $\beta_i(z)=z$, $\delta_i(z)=w^{i-1}z$,  $i=1,\cdots,N$, and $a(z)=f(z)$, $b(z)=0$.  Obviously the corollary can be used since $\beta_i^{-1}(z)=z$, $\delta_i^{-1}(z)=w^{N+1-i}z$, and 
\[
(\beta_i^{-1}\circ \delta_i-\beta_j^{-1}\circ \delta_j)(z)=(\delta_i-\delta_j)(z)=(w^{i-1}-w^{j-1})z
\]
is an automosphim of $\mathbb{C}$ for all $i\neq j$. 
{\hfill $\Box$}

Another special case of Wilson's equation is 
\begin{equation}\label{Ski-ordinary}
\sum_{i=1}^{m}f_i(b_ix+c_iy)= \sum_{i=1}^m f_i(b_ix) +   \sum_{i=1}^m f_i(c_iy),
\end{equation}
which is a linearized form of the
Skitovich-Darmois functional equation
\begin{equation*}
\prod_{i=1}^m \widehat{\mu_i}(b_ix+c_iy) = \prod_{i=1}^m \widehat{\mu_i}(b_ix)  \prod_{i=1}^m  \widehat{\mu_i}(c_iy),
\end{equation*}
an equation which is connected to the characterization  problem of Gaussian distributions (see, for example, [Linnik \cite{Linn}, Ghurye-Olkin \cite{Gu_O}, \cite{KLR}]):

\begin{corollary}[Linearized Skitovich–Darmois functional equation]  \label{SD} Assume that  $(S,+)$ and $(H,+)$  are commutative groups 
and let $\beta_i,\delta_i:S\to S$ be automorphisms such that $\beta_i^{-1}\circ \delta_i-\beta_j^{-1}\circ \delta_j$ is invertible whenever $i\neq j$. If the functions $f_i:S\to H$, $i=1,\cdots, n$ solve the functional equation
\begin{equation}\label{LSDFE}
\sum_{i=1}^nf_i(\beta_i(x)+\delta_i(y))=\sum_{i=1}^nf_i(\beta_i(x))+\sum_{i=1}^nf_i(\delta_i(y)),
\end{equation}
then $P(x)=\sum_{i=1}^nf_i(\beta_i(x))$ and $Q(y)= \sum_{i=1}^nf_i(\delta_i(y))$ are polynomial functions of degree $\leq n$. 
\end{corollary}

\noindent  \textbf{Proof. } Use Corollary \ref{Wilson} with $a(x)=P(x)$ and $b(y)=Q(y)$.

{\hfill $\Box$}

Equation \eqref{LSDFE} has been studied in  great detail by Feld'man \cite{Feld} for functions defined on locally compact commutative groups.

\begin{remark} As a final remark we would like to mention that although we have formulated all results in this paper for ordinary functions defined on commutative groups or semigroups, the same proofs can be translated to the distributional setting. In particular, ordinary polynomials of degree $\leq m$ are the unique solutions of Aichinger's equation \eqref{aichinger} when $f\in \mathcal{D}(\mathbb{R}^d)'$ (so that $f(x_1+\cdots+x_{m+1}) \in \mathcal{D}(\mathbb{R}^d\times \cdots ^{m+1 \text{times}}\times \mathbb{R}^d)'$) and $g_i(x_1,\cdots,\widehat{x_i},\cdots,x_{m+1}) \in  \mathcal{D}(\mathbb{R}^d\times \cdots ^{m+1 \text{times}}\times \mathbb{R}^d)'$ but does not depend on $x_i$). The proof of this result is a direct consequence of the fact that the translation and difference operators are well defined for distributions and inherit in the distributional framework the properties they have when applied to ordinary functions, and Fréchet's functional equation also characterizes polynomials in a distributional sense \cite{A_NFAO}. 
\end{remark}

\section{Declarations}

 \noindent \textbf{Ethical Approval}
 
 Not applicable
 \vspace{0.5cm} 
 
 \noindent \textbf{Competing interests }
 
The author has no relevant financial or non-financial interests to disclose.

\vspace{0.5cm} 
 
 \noindent \textbf{Authors' contributions }
 
Not applicable \vspace{0.5cm} 
 
 \noindent \textbf{Funding }
 
Not applicable \vspace{0.5cm} 
 
 \noindent \textbf{Availability of data and materials }

Not applicable


\begin{thebibliography}{1}
\bibitem{AM} {\sc E. Aichinger and J. Moosbauer, } Chevalley-Warning type results on abelian groups. J. of Algebra 569 (2021) 30-66.
%
%\bibitem{AA_AM} {\sc A. Aksoy, J. M. Almira, } On Montel and Montel-Popoviciu Theorems in several variables, Aequationes Mathematicae, \textbf{89} (5) (2015) 1335--1357.

\bibitem{almira}  J. M. Almira, Aichinger equation on commutative semigroups, \textit{arXiv:2201.07797 }, 2022.

\bibitem{almirajmaa} {\sc J. M. Almira, } Characterization of polynomials as solutions of certain functional equations, J. Math. Anal. Appl. 459 (2018) 1016–1028

\bibitem{A_NFAO} {\sc J. M. Almira, }  Montel's theorem and subspaces of distributions which are $\Delta^m$-invariant, Numer. Functional Anal. Optimiz. \textbf{35} (4) (2014) 389--403.

%\bibitem{A_popo} {\sc J. M. Almira, } On Popoviciu-Ionescu functional equation, Annales Mathematicae Silesianae,  \textbf{30} (2016) 5--15.

\bibitem{A_JJA} {\sc J. M. Almira, } On Loewner's characterization of polynomials, Ja\'{e}n Journal on Approximation,  \textbf{8} (2)  (2016) 175--181.
%
%\bibitem{AK_CJM} {\sc J. M. Almira, K. F. Abu-Helaiel, } On Montel's theorem in several variables, Carpathian Journal of Mathematics,  \textbf{31} (2015), 1--10.

%\bibitem{AK_ATPSFE} {\sc J. M. Almira, K. F. Abu-Helaiel, }  A note on invariant subspaces and the solution of some classical functional equations, Annals of the Tiberiu Popoviciu Seminar of Functional Equations, Approximation and Convexity \textbf{11} (2013) 3-17.
%
%\bibitem{AS_LCFE}  {\sc J. M. Almira, E. V. Shulman, } On certain generalizations of the Levi-Civita and Wilson functional equations,  Aequationes Mathematicae, Aequationes Math. \textbf{91} (5) (2017) 921–931.
%
%\bibitem{AS_AM} {\sc J. M. Almira, L. Sz\'{e}kelyhidi, } Local polynomials and the Montel theorem, Aequationes Mathematicae, \textbf{89} (2015) 329-338. 

%\bibitem{AS_DM} {\sc J. M. Almira, L. Sz\'{e}kelyhidi, } Montel--type theorems for exponential polynomials,  Demonstratio Mathematica, \textbf{49} (2) %(2016) 197-212.


%\bibitem{anselone} {\sc P. M. Anselone, J. Korevaar, } Translation invariant subspaces of finite dimension, \emph{Proc. Amer. Math. Soc. } \textbf{15} (1964), 747-752.
%
%\bibitem{Cra} {\sc H. Cram\'{e}r, } \emph{Random variables and probability distributions, } Cambridge Tracts in Maths. and Math. Physics \textbf{36} Cambridge University Press, London, 1937. 

%\bibitem{E} {\sc M. Engert, } Finite dimensional translation invariant subspaces, Pacific J. Maths. \textbf{32} (2) (1970) 333-343.

%\bibitem{fe} {\sc I. Feny\"{o}, } Sur les \'{e}quations distributionnelles, in \emph{Functional Equations and Inequalities,} Ed. B. Forte, C.I.M.E. summer %schools, \textbf{54}, Springer, 1971 (Reprinted in 2010), 45-109.

\bibitem{Feld} {\sc G. Feldman, } \emph{Functional equations and characterization problems on locally compact Abelian groups, } EMS, 2008.

\bibitem{Gu_O} {\sc S. G. Ghurye, I. Olkin, } A characterization of the multivariate normal distribution,  Ann. Math. Statist. \textbf{33} (2) (1962) 533--541. 

%\bibitem{HW} {\sc G. H. Hardy, E. M. Wright, } \emph{An Introduction to the Theory of Numbers. Fifth edition.} The Clarendon Press, Oxford University Press, New York, 1979.

\bibitem{h1} {\sc S. Haruki, } On the mean value property of harmonic and complex polynomials, Proc. Japan Acad. Ser. A, \textbf{57}  (1981) 216-218. 

\bibitem{h2} {\sc S. Haruki, } On the theorem of S. Kakutani-M. Nagumo and J.L. Walsh for the mean value property of harmonic and complex polynomials, Pacific J. Math. \textbf{94} (1) (1981) 113-123.

\bibitem{h3} {\sc S. Haruki, } On two functional equations connected with a mean-value property of polynomials, Aequationes Math. \textbf{6} (1971) 275-277.

\bibitem{KLR} {\sc A.M. Kagan, Yu. V. Linnik, C. R. Rao, } \emph{Characterization problems in Mathematical Statistics, } John Wiley \& Sons, 1973.

\bibitem{kn} {\sc S. Kakutani, M. Nagumo, } About the functional equation $\sum_{v=0}^{n-1}f(z+e^{(\frac{2 v \pi}{2})i})=nf(z)$, Zenkoku Shij\^{o} Danwakai, \textbf{66} (1935) 10-12 (in Japanese).

%\bibitem{leland} {\sc K. O. Leland, } Finite dimensional translation invariant spaces, \emph{Amer. Math. Monthly} \textbf{75} (1968) 757-758.


\bibitem{Linn} {\sc Yu. V. Linnik,} \emph{A remark on the Cramer's theorem on the decomposition of the normal law}, Theory Prob. Appl. \textbf{1} (1956), no. 4, 435--436. 

%\bibitem{Lo} {\sc C. Loewner, } On some transformation semigroups invariant under Euclidean and non-Euclidean isometries, J. Math. Mech.  8 (1959) 393-409. 

\bibitem{MaPe} {\sc A. M. Mathai, G. Pederzoli, } \emph{Characterizations of the Normal Probability Law, } Wiley Eastern Limited, 1977.

%\bibitem{S1} {\sc E. Shulman, } Subadditive set-functions on semigroups, applications to group representations and functional equations, \emph{J. Funct. %Anal.} \textbf{263} (2012) 1468-1484.

%\bibitem{S} {\sc E. Shulman, } Decomposable functions and representations of topological semigroups, \emph{Aequat. Math.} \textbf{79} (2010) 13-21.

%\bibitem{S_PhD} {\sc E. Shulman, }  \emph{Functional equations of homological type,}  PhD Thesis, Moscow State Pedagogical University (1994).


\bibitem{Sz1} {\sc L. Szekelyhidi} \emph{Convolution type functional equations on Topological Abelian Groups, } World Scientific, 1991.

\bibitem{Sz} {\sc L. Szekelyhidi} \emph{Discrete spectral synthesis and its applications, } Springer, 2006.

%\bibitem{V} {\sc V. S. Vladimirov, } \emph{Generalized Functions in Mathematical Physics,} 1979.

%\bibitem{W} {\sc M. Waldschmidt, } \emph{Topologie des Points Rationnels, } Cours de Troisi\`{e}me Cycle 1994/95 Universit\'{e} P. et M. Curie (Paris VI), %1995.

\bibitem{wal} {\sc J. L. Walsh, } A mean value theorem for polynomials and harmonic polynomials, Bull. Amer. Math. Soc. \textbf{42} (1936) 923-930. 

\bibitem{wilson} {\sc W.H. Wilson, }  On a certain general class of functional equations. Am. J. Math. \textbf{40} (3) (1918) 263–282. 


\end{thebibliography}
\end{document}